\documentclass[12pt]{amsart}
\usepackage{amsmath,amssymb}
\usepackage{color}
\theoremstyle{plain}
\newtheorem{thm}{Theorem}[section]
\newtheorem{lem}[thm]{Lemma}

\newtheorem{cor}[thm]{Corollary}
\newtheorem{pro}[thm]{Proposition}
\newtheorem{ex}[thm]{Example}
\newtheorem{que}[thm]{Question}
\newtheorem{prob}[thm]{Problem}

\newtheorem{df}[thm]{Definition}

\newcommand{\ind}{\mbox{{\rm ind}\,}}

\newcommand{\Int}{\mbox{{\rm Int}}}
\newcommand{\Cl}{\mbox{{\rm Cl}}}

\begin{document}
\title{On (strongly) ($\Theta$-)discrete homogeneous spaces}

\author{Vitalij A.~Chatyrko and Alexandre Karassev}

\begin{abstract} We introduce the classes of (strongly) ($\Theta$-)discrete homogeneous spaces. We discuss the relationships of these classes to other classes of spaces possessing homogeneity-related properties, such as (strongly) ($n$-)homogeneous spaces. Many examples are given distinguishing discrete homogeneity and other types of homogeneity.
\end{abstract}

\makeatletter
\@namedef{subjclassname@2020}{\textup{2020} Mathematics Subject Classification}
\makeatother



\maketitle

\section{Introduction} 

A family of subsets $\mathcal D$ of a topological space $X$ is 
{\it discrete} if each point $x \in X$ has a neighborhood  $Ox$ which  intersects at most one  set from $\mathcal D$. If each element of $\mathcal D$ is a singleton, we obtain the definition of a discrete subset of $X$. Note that if $X$ is $T_1$ then each discrete subset of $X$ is closed.

Recall also \cite{BBHS} that a subset $D$ of a  topological space $X$ is 
{\it $\Theta$-discrete} if each point $x \in X$ has a neighbourhood  $Ox \subset X$ whose closure
$\Cl Ox$ contains at most one point of the set $D$.

It is clear that each $\Theta$-discrete  subset of a topological space  is discrete. There exist Hausdorff spaces in which some discrete subsets are not $\Theta$-discrete (\cite[Example 1]{BBHS}). However, in regular Hausdorff spaces  $\Theta$-discrete  subsets and  discrete subsets are the same.

In \cite{BBHS} the authors proved the following interesting result.

\begin{thm}{\rm (}\cite[Theorem 1 and Corollary 2]{BBHS}{\rm )}  Let $X$ be 
the connected countable Hausdorff space $B$ due to Bing (\cite{B}, see also the appendix for the definition) or
 the connected, locally connected, countable Hausdorff space $R$ due to Ritter (\cite{R}, see  Section 3 for the definition).
 Then any bijection $f\colon A \to B$ between $\Theta$-discrete subsets $A$ and $B$ of $X$ extends to a homeomorphism  $f\colon X \to X$.
\end{thm}

This result and the definition of (strongly) n-homogeneous spaces (see further in the text) motivates the following definition (for a set $S$ by $|S|$ we denote its cardinality).

\begin{df}\label{main_def}  A Hausdorff topological space $X$ is called 
\begin{itemize}
\item {\it discrete homogeneous} ($DH$), if for any two discrete subsets $A$ and $B$ of $X$ with $|A|=|B|$ there exists a homeomorphism 
$f$ of $X$ onto itself such that $f(A)= B$;

\item {\it strongly discrete homogeneous} ($sDH$), if for any two discrete subsets $A$ and $B$ of $X$ and any bijection $f\colon A\to B$ this bijection extends to a homeomorphism of the whole $X$ onto itself.

\end{itemize}

\end{df}

Replacing in the above definition discrete subsets with $\Theta$-discrete, we obtain the definitions of $\Theta$-discrete homogeneous ($\Theta$-DH), and strongly $\Theta$-discrete homogeneous (s$\Theta$-DH) spaces, respectively.

Note that the Bing space $B$ and the Ritter space $R$ are examples of  s$\Theta$-DH spaces. It is also easy to see that each s($\Theta$-)DH space is ($\Theta$-)DH, and each (s)DH space is (s)$\Theta$-DH. Note, however, that the Bing space $B$ is not a DH space (\cite[Example 2]{BBHS}). It is natural to ask the following question.

 \begin{que} 
Is the space $R$ an (s)DH-space?

 \end{que}

Recall (\cite{Bu}) that a topological space $X$ is called {\it $n$-homogeneous} for a positive integer $n$ if
for any two  subsets $A$ and $B$ of $X$ of cardinality $n$ there exists a homeomorphism $f\colon X \to X$, such that $f(A)=B$. The $1$-homogeneous spaces are called homogeneous. Further, $X$ is called {\it strongly $n$-homogeneous} if
any bijection  $f\colon A \to B$ between  subsets $A$ and $B$ of $X$ of cardinality $n$ extends to a homeomorphism  $f\colon X \to X$.

Since any finite subset of a Hausdorff space is $\Theta$-discrete we have the following fact.

\begin{pro}(\cite[Corollary 1 for the space $B$ and Corollary 3  for the space $R$]{BBHS}) Each (s)$\Theta$-DH space (resp. (s)DH space) is (strongly) $n$-homogeneous for every positive integer $n$. 
\end{pro} 
 
In this article we will discuss the introduced  classes of topological spaces, and their relationships to other classes of spaces with homogeneity-type properties. We prove that several classical spaces possess the (strong) discrete homogeneity.

 The reader can find all standard notions mentioned in the text in \cite{E} or \cite{M}. Everywhere below, $\mathbb N$, $\mathbb Q$, $\mathbb P$, $\mathbb C$, $\mathbb R$, $S^1$, $\mathbb S$ denote the sets of natural, rational, and irrational numbers, the Cantor set, real line, the circle, and the Sorgenfrey line, respectively.

\section{The case of regular Hausdorff spaces}

The following statement is evident.
\begin{pro} Let $X$ be a regular Hausdorff space. Then $X$ is an  (s)$\Theta$-DH space iff it is  an  (s)DH space. 
\end{pro}

It is easy to see that the real line $\mathbb R$ is not strongly 3-homogeneous (although it is strongly 2-homogeneous) and the circle $S^1$ is not strongly 4-homogeneous (although it is strongly 3-homogeneous). Therefore both these spaces are not
sDH spaces. Note also that $S^1$ is an example of a DH space that is not sDH.

\begin{pro} Any infinite discrete space $X$ is strongly $n$-homogeneous for every positive integer $n$, but not a  DH space.
\end{pro} 

\begin{proof} It is clear that  $X$ is strongly $n$-homogeneous for every positive integer $n$. We will show that $X$ is not a DH space. Let $a, b$ be two different points of $X$. Put $A = X \setminus \{a\}$ and $B = X \setminus \{a, b\}$. Clearly, $|A|=|B|$, and $A, B$ are discrete subsets of $X$.  However, there is no bijection $f\colon X\to X$ such that $f(A)=B$.
 
\end{proof}

A Hausdorff space is called {\it strongly homogeneous} (see \cite{vM}) provided that all nonempty clopen subspaces are homeomorphic. 

The spaces $\mathbb C$, $\mathbb Q$,  $\mathbb P$, ${\mathbb C}\times\mathbb Q$, ${\mathbb C}\times\mathbb P$ and $\mathbb S$ are strongly homogeneous (for the product  ${\mathbb C}\times\mathbb Q$ see \cite{AU}, and for the product  ${\mathbb C}\times\mathbb P$ see \cite{vM}).

It is easy to see that each zero-dimensional in the sense of dimension $\ind$ first countable strongly homogeneous space is homogeneous. See \cite[p. 383]{E} for the definition of the small inductive dimension $\ind$.

The following interesting result was obtained in \cite{vM}.

\begin{thm}(\cite[Theorem 3.1]{vM})\label{vM_theorem}   Let $X$ be 
a strongly homogeneous separable metrizable zero-dimensional space $X$ and let $F$ and $G$
be closed nowhere dense subsets of $X$. 
If $f\colon F \to G$ is a homeomorphism then $f$ can be extended to an automorphism   $\bar{f}$ of $X$.
\end{thm}

From Theorem \ref{vM_theorem} we immediately obtain the following.

\begin{pro}\label{the Cantor set} Every strongly homogeneous separable metrizable zero-dimensional space is sDH.
In particular, the spaces  $\mathbb C$, $\mathbb Q$,  $\mathbb P$, ${\mathbb C}\times\mathbb Q$, ${\mathbb C}\times\mathbb P$ are sDH.
\end{pro}

In \cite{vM} it was observed that the zero-dimensional separable metrizable topological group
$H$ from \cite{vD} (in fact, $H$ is a subgroup of the circle group) is not strongly homogeneous. 
The necessary property of $H$ for the above fact is the following:

{\it if $F$ and $G$ are homeomorphic closed subspaces of $H$ then $H \setminus F$ and $H \setminus G$ are homeomorphic (note that we do not require that $F$ and $G$ are proper closed subsets of $H$).}

We will prove that $H$ is sDH.

Let us say that a Hausdorff space $X$ possesses {\it a property $\mathcal D$} if for any proper clopen homeomorphic subspaces $F$ and $G$ of $X$ the subspaces $X \setminus F$ and $X \setminus G$ are also homeomorphic.

It is easy to see that the group $H$ as well as any strongly homogeneous space possesses
 property $\mathcal D$.

 Recall further that a  space $X$ is called {\it collectionwise normal} if $X$ is $T_1$ and for every discrete family $\{F_s\}_{s \in S}$ of closed subsets of $X$ there exists a discrete family $\{V_s\}_{s \in S}$  of open subsets of $X$ such that $F_s \subset V_s$ for every $s \in S$. It is a well-known fact that every paracompact space is collectionwise normal (\cite[p. 305]{E}). 

\begin{lem}\label{coll_norm} Let $X$ be a collectionwise normal space with $\ind X =0$, and $A$ be a discrete subset of $X$. 
Then there exists disjoint family $\{V_a\colon a \in A\}$ of clopen
subsets of $X$ such that $a \in V_a$ for each $a \in A$, and 
the set $X \setminus \cup_{a \in A} V_a$ is also clopen. 
\end{lem}

\begin{proof} Since the set $A$ is discrete and the space $X$ is $T_1$,  the family $F_a = \{ \{a\}\colon a \in A \},$ is discrete and consists of closed subsets of $X$. Therefore there exists a discrete family $\{V_a\}_{a \in A}$  of open subsets of $X$ such that $a \in V_a$ for every $a \in A$. We may assume  that all $V_a, a \in A,$ are clopen because $\ind X = 0$.
Since  the family $\{V_a\}_{a \in A}$ is discrete, the set
$\cup_{a \in A} V_a$ is closed. Hence  the set $X \setminus \cup_{a \in A} V_a$ is also clopen. 
\end{proof}

\begin{thm}  Let $X$ be a collectionwise normal homogeneous  space with  $\ind X =0$ and with property $\mathcal D$. Then $X$ is sDH.
\end{thm}

\begin{proof}  Consider two discrete
subsets $A, B$ of $X$ and any bijection $f\colon A \to B$. Since the space $X$ is  collectionwise normal and $\ind X=0$
we can apply Lemma \ref{coll_norm} to each of the sets $A$ and $B$. So 
\begin{itemize}
\item[(a)] there exists a disjoint family $\{V_a\colon a \in A\}$ of clopen
subsets of $X$ such that $a \in V_a$ for each $a \in A$, and
the set $X \setminus \cup_{a \in A} V_a$ is proper clopen, and 
\item[(b)] there exists a disjoint family $\{W_b\colon b \in B\}$ of clopen
subsets of $X$ such that $b \in W_b$ for each $b \in B$, and 
the set $X \setminus \cup_{b \in B} W_b$ is proper clopen.
\end{itemize}

Furthermore, since $X$ is homogeneous we can assume that for each $a \in A$ 
the subspaces $V_a$ and $W_{f(a)}$ are homeomorphic and for a homeomorphism
$f_a\colon V_a \to  W_{f(a)}$ we have
$f_a(a) = f(a) = b$. Set now $F = \cup_{a \in A} V_a$ and $G = \cup_{b \in B} W_b$.
Note that $F$ and $G$ are clopen, proper and homeomorphic. So by the property $\mathcal D$
the subspaces $X \setminus F$ and $X \setminus G$ are also homeomorphic.
This easily implies that there exists an automorphism  $\bar{f}$ of $X$ which extends $f$.
\end{proof}

 \begin{cor} 
Let $X$ be a collectionwise normal first countable strongly homogeneous topological space with $\ind X =0.$ Then $X$ is sDH. 

\end{cor}

 \begin{cor} 
The Sorgenfrey line $\mathbb S$ is sDH.
\end{cor}

 \begin{cor} The topological group $H$ from \cite{vD} is sDH.
\end{cor}

\begin{que} 
\begin{itemize}
\item[(a)] Is $H^n$ sDH for each positive integer $n$?
\item[(b)] Is $\mathbb S^n$ sDH for each positive integer $n$?
\end{itemize}
\end{que}

Theorem \ref{vM_theorem} 
 motivates the following definition.

\begin{df}\label{main_def}  A separable metrizable zero-dimensional topological space $X$ is called 
 {\it nowhere dense homogeneous} (NDH), if for any two closed nowhere dense subsets $F$ and $G$ of $X$ and any homeomorphism  $f\colon F \to G$ there exists an automorphism $\bar{f}$ of $X$ which extends $f$.
\end{df}

The discussion above indicates that each  separable metrizable zero-dimensional strongly homogeneous  space $X$ is NDH, and each NDH space is sDH.

\begin{que}
\begin{itemize}
\item[(a)] Is there NDH space which is not strongly homogeneous?
\item[(b)] Is there separable metrizable zero-dimensional sDH space which is not NDH?
\end{itemize}
\end{que}

Since the group $H$ is sDH but not strongly homogeneous, $H$ provides a positive answer to one of the questions.

Let us continue with the case $\ind \geq 1.$ 
 
 For  $x$ in a metric space $X$ and $r>0$ let $B(x,r)$ be the open ball of radius $r$ centred at $x$. We will also let $B(x,0) = \{x\}.$

 Next, we will prove that $\mathbb R^n$ is sDH for all $n\ge 2$. Since the proofs differ for $n=2$ and $n>2$, we divide the proofs into two propositions below. Note that these results show that the product of two DH spaces (in this case, $\mathbb R$, which is not sDH) can be sDH.

  \begin{pro} $\mathbb R^2$ is 
 an  sDH space.
\end{pro}

\begin{proof} First, we will show that $\mathbb R^2$ is DH. Let $A$ be a countable discrete subset of $\mathbb R^2$. We will construct a homeomorphism $h\colon \mathbb R^2 \to \mathbb R^2$ such that $f(A) = \mathbb N$, where $\mathbb N$ is viewed as a subset of $\mathbb R^2$. Using general position argument, we may assume that each circle centred at the origin contains at most one point from $A$. Hence, $A = \{ a_k\mid k\ge 1\}$ and we may assume that $0 < \|a_1\| < \|a_2\| < \|a_3\| < \dots$ Let $\epsilon _k =  \min\{\|a_{k+1}\| - \| a_k\|, \|a_k\| - \| a_{k-1}\|\}$, $k\ge 1$ (here, we let $a_0 =(0,0)$). Let $U_k$ be  $\epsilon _k /3$ - neighbourhood of a circle of radius $\| a_k\|$ centred at the origin. For each $k$, there exists a homeomorphism $f_k$ of $\mathbb R^2$ moving $a_k$ to $(\|a_k\|,0)$  and being identity outside $U_k$. We can now define the limit $f$ of $f_k\circ f_{k-1} ... \circ f_1$. Clearly, $f$ is a homeomorphism. Note that $A' = f(A)$ is a discrete subset of the positive part of first coordinate axis in $\mathbb R^2$. Clearly, there exists a homeomorphism $g$ of $\mathbb R^2$ such that $g(A') = \mathbb N$. Now we let $h = g\circ f$. 

Further, let $A$ and $B$ be two discrete countable subsets of $\mathbb R^2$. Again, using general position argument, we may assume that $A\cap B = \varnothing$. We want to extend a bijection between $A$ and $B$ to a homeomorphism of $\mathbb R^2$. For this, we first consider a homeomorphism  
of $\mathbb R^2$ that maps $A\cup B$ onto $\mathbb N$ (it exists because of the first part of the proof). Thus, we reduce the task of extending bijection between $A$ and $B$ to the following: for two disjoint infinite subsets $C$ and $D$ of $\mathbb N$ such that $C\cup D = \mathbb N$ and a bijection of $C$ onto $D$, find a homeomorphism of $\mathbb R^2$ extending this bijection. So let $C=\{ c_k\mid k\ge 1\}$, $D=\{d_k\mid k\ge 1\}$, and we want to extend a bijection sending $c_k$ to $d_k$ for each $k$. Fix a countable collection of rays $R_k$, $k\ge 1$, in $\mathbb R^2$, such that $R_i\cap R_j = (0,0)$ for all $i\ne j$. Further, for each $k$, let $U_k$ and $V_k$ be $1/3$ neighbourhoods of circles of radii $\|c_k\|$ and $\| d_k\|$, respectively, centred at the origin. There exists a homeomorphism $f_k$, $k=1,2,\dots$, such that $c'_k = f_k(c_k)\in R_k$, $d'_k = f_k(d_k)\in R_k$, $\|c'_k\| = \|c_k\|$, $\|d'_k\|= \| d_k\|$ and $f_k$ is the identity outside of $U_k\cup V_k.$ Let  $f$ be the limit of $f_k\circ f_{k-1} ... \circ f_1$. Clearly, $f$ is a homeomorphism. Note that the collection of  
    straight line segments $[c'_k, d'_k],$ joining $c'_k$ and $d'_k$, is discrete. For each $k$, let $g_k$ be a homeomorphism of $\mathbb R^2$ such that $g_k (c'_k) = d'_k,$ which is the identity outside of some ``narrow" neighbourhood of  $[c'_k, d'_k].$ We can define a homeomorphism $g$ which is the limit of $g_k\circ g_{k-1} ... \circ g_1$. The homeomorphism $f^{-1}\circ g\circ f$ is the one we are looking for.
\end{proof}

 \begin{pro} Every finite-dimensional Euclidean space $\mathbb R^n, n >2,$ is 
 an  sDH space.
\end{pro}

\begin{proof} Indeed, let $A$ and $B$ be two discrete subsets of $\mathbb R^n$ of the same cardinality. Hence both $A$ and $B$ are either finite or countable. We will assume $A$ and $B$ are countable (the finite case is similar). Therefore $A = \{a_k\colon k \geq 1\}$ and $B = \{b_k\colon k \geq 1\}$.  We want to extend a bijection sending $a_k$ to $b_k$ for each $k \geq 1$.  For  any $k$ there exists a path  $P_k$ joining $a_k$ and $b_k$, that avoids the ball of radius $\frac{1}{2}\min\{\|a_k\|, \|b_k\|\}$ centred at the origin, and also avoids all other $a$'s and $b$'s. Using general position argument in $\mathbb R^n$, $n>2$, we may also assume that $P_i\cap P_j = \varnothing$ for $i\ne j$. There exists a homeomorphism $f_k\colon \mathbb R^n\to \mathbb R^n$ that maps $a_k$ to $b_k$ and is the identity outside some ``narrow" neighbourhood of the path $P_k$. We can now define the limit of $f_k\circ f_{k-1} ... \circ f_1$, which is the homeomorphism extending the bijection between $A$ and $B$.\end{proof}

\begin{pro} The separable Hilbert space $l_2$ is sDH.
\end{pro}

\begin{proof} As in the proof of the previous propositon, consider two countable subsets $A = \{a_k\colon k \geq 1\}$ and $B = \{b_k\colon k \geq 1\}$ of $l_2$. We  want to extend a bijection sending $a_k$ to $b_k$ for each $k \geq 1.$
 
  For each $x\in l_2$ there exists $r_x>0$ such that $B(x, r_x) \cap (A\cup B) =\{ x\} \cap (A\cup B)$. Additionally, we may assume that $r_{a_k} < 1/k$ and $r_{b_k}<1/k$ for all $k$.  Let $\mathcal U =\{ B(x, r_x) \mid x\in l_2\}$. 
 
 Join each  $a_k$ to $b_k$ by a straight line segment. This collection of segments can be viewed as a map $g\colon \oplus_{k=1}^\infty ([0,1]\times \{k\})\to l_2$. Therefore there exists a map $g'\colon \oplus_{k=1}^\infty ([0,1]\times \{k\})\to l_2$ such that the family $\{ g'([0,1]\times \{k\})\mid k\in \mathbb N\}$ is discrete. Each $g'([0,1]\times \{ k\})$ is a path such that $g'(0)\in B(a_k,r_{a_k})$ and $g'(1)\in B(b_k, r_{b_k})$. By adding straight line segments from $a_k$ to $g'(0)$ and from $b_k$ to $g'(1)$ we obtain a discrete family of paths $P_k$ in $l_2$, joining $a_k$ and $b_k$. Moreover, we may even assume that each $P_k$ is piecewise linear arc. 
There exists a homeomorphism $f_k\colon  l_2\to l_2$ that maps $a_k$ to $b_k$ and is identity outside some ``narrow" neighbourhood $V_k$ of the path $P_k$, such that the family $\{ V_k\colon k\ge 1\}$ is discrete. As in the proof of the previous proposition, the limit of $f_k\circ f_{k-1} ... \circ f_1$ is the desired homeomorphism.
\end{proof}

Recall that a space $X$ is called strongly locally homogeneous (SLH) if it has a basis of open sets $\mathcal B$ such that for every $U\in\mathcal B$ and any two points $x$ and $y$ in $U$ there exists a homeomorphism $f\colon X\to X$ such that $f(x) =y$ and $f$ is identity on $X\setminus U$. Note that any $n$-manifold for $n>1$ and any Hilbert space manifold are examples of SLH spaces. It was shown in \cite{Ba} that  a connected SLH space, no two-point subset of which disconnects it, is strongly $n$-homogeneous for any $n$.

\begin{ex}\label{nondh} Let $M = \mathbb R\times S^1.$ Then $M$ is a $2$-manifold which is strongly $n$-homogeneous for any $n$ and SLH but not DH
\end{ex}
\begin{proof} Being a $2$-manifold, $M$ is both strongly $n$-homogeneous and SLH. Fix a point $s\in S^1$ and let $A = \{ (n,s) \colon n\in \mathbb N\}$ and  $B = \{ (\pm n,s) \colon n\in \mathbb N\}.$  We claim that there is no autohomeomorphism of $M$ that maps $A$ onto $B$. Indeed, suppose $f\colon M\to M$ is such homeomorphism.  Let $N = [0,\infty)\times S^1.$ Note that $f(N)$ is a submanifold of $M$, containing $B$, with the boundary being simple closed curve $C = f(\{0\}\times S^1)$. If $C$ is nullhomotopic, then $C$ is a common boundary of a topological disk and a manifold, homeomorphic to a sphere with three holes. Thus, in this case $f(N)$ cannot be homeomorphic to $N$. Otherwise, $C$ divides $M$ into two submonifolds, none of which can contain the whole $B$, contradiction. 
\end{proof}

 \begin{que} Which metrizable manifolds are DH or sDH? 
 \end{que}

 \begin{prob} Find conditions for SLH space to be (s)DH. 
 \end{prob}
 
 In this question, one can impose various additional conditions such as being completely metrizable and separable, connected, arc-wise connected, locally compact, etc.

Since every discrete subspace of a compact Hausdorff space is finite we have the following
statement.

\begin{pro} Let $X$ be a compact Hausdorff space. Then $X$ is a (s)DH space iff it is (strongly) $n$-homogeneous for every positive integer $n$. 
\end{pro}

Besides the Cantor set $\mathbb C$ (see Proposition \ref{the Cantor set}), an obvious  example of a zero-dimensional sDH compact Hausdorff space is the Alexandroff double arrow space \cite{AU2}. Note also, that since $\mathbb C$ is the only homogeneous metrizable  zero-dimensional infinite compactum, it is also the only metrizable  zero-dimensional infinite compactum which is sDH. In general case, however, the following problem seems interesting.

\begin{prob} Characterize all zero-dimensional sDH compacta.
\end{prob}

Recall that a closed subset $A$ of a metric space $X$ is called a $Z$-set if for any $\epsilon >0$ there exists a map $f\colon X \to X$ such that $f$ is $\epsilon$-close to the identity map and $f(X)\cap A = \varnothing$. Let $\mu^k$ denote the Menger cube for each $k=1,2,\dots$, and $Q=[0,1]^\omega$ be the Hilbert cube.  It is well-known that any finite subset of $\mu^k$ or $Q$ is a $Z$-set in the respective space, and any homeomorphism between two $Z$-sets of $\mu^k$ or $Q$ can be extended to a homeomorphism of the whole space (see, e.g. \cite{Ch} for the case of the Hilbert cube, and \cite{A, Be} for the case of Menger compacta). Consequently, $Q$ and $\mu^k$ are examples of sDH compacta.

Moreover, it is well-known \cite{A} that $\mu^1$ and $S^1$ are the only homogeneous $1$-dimensional Peano continua (by a continuum we mean metrizable compact connected space). By the result of \cite{U}, any $2$-homogeneous continuum is locally connected. Therefore we have the following observation.

\begin{pro}
$\mu^1$ and $S^1$ are the only DH $1$-dimensional continua.
\end{pro}

Since as was observed earlier $S^1$ is not sDH but $\mu^1$  is sDH, we get immediately the following corollary.

\begin{cor} $\mu^1$ is the only sDH $1$-dimensional continuum
\end{cor}

Note that  according to \cite{KKT}, the product $\mu ^1 \times \mu ^1$ is not $2$-homogeneous, and hence not DH. Moreover, it is easy to see that $\mathbb R^n\times \mathbb Q$, $n \geq 1,$ is also not $2$-homogeneous.
Therefore the property (s)DH is not preserved by products.

In view of the above discussion about $n$-homogeneity, as well as  Example \ref{nondh}, the following question appears to be interesting.
 
 \begin{que} Let $X$ be a connected separable metrizable (complete, locally compact, locally connected etc.)  space, that is (strongly) $n$-homogeneous for any $n.$ What additional conditions on $X$ imply that $X$ is (s)DH?
 \end{que}

\section{The case of non-regular Hausdorff spaces} 
Since every countable regular Hausdorff space cannot be connected, the Bing space $B$ and the Ritter spaces $R$ are non-regular.

Let $X$ be a topological space with the topology $\tau$. 

Recall (cf. \cite{MRV}) that the topology 
$\hat{\tau} \subseteq \tau$ on the set $X$ generated by the base consisting of all regular open subsets
of the topological space $X$ is called the semi-regularization of $\tau$.
(An open set $O$ of a topological space $X$ is called  regular open if $O$ coincides with the interior of its  closure $\Cl (O)$ in $X$.) 

We will denote the space $(X, \hat{\tau})$ by $\hat{X}$. Recall (\cite{R}) that 
$\hat{B} = R$. 

Below we recall some facts about  the semi-regularization topologies.

\begin{pro}(cf. \cite{MRV}) Let $X$ be a topological space with the topology $\tau$. Then we have the following.
\begin{itemize}

\item[(a)] $\hat{X}$ is Hausdorff iff $X$ is Hausdorff;

\item[(b)] $\hat{X}$ is connected iff $X$ is connected;

\item[(c)] if $X$ is locally connected then $\hat{X}$ is also locally connected but the converse is false.

\end{itemize}

\end{pro}

Note that the Bing space $B$ is not locally connected.

The following lemma is evident.

\begin{lem}\label{thetasemiregularization} Let $X$ be a Hausdorff topological space.  
If $D$ is a  $\Theta$-discrete (resp. discrete) subset of $\hat{X}$ then $D$ is a $\Theta$-discrete (resp. discrete)
subset of $X$.
\end{lem}

\begin{thm}(For the spaces $B$ and $R$, \cite{BBHS}) If $X$ is an s$\Theta$-DH space then $\hat{X}$ is also an s$\Theta$-DH space.
\end{thm}

\begin{proof} Let $A$ and $B$ be $\Theta$-discrete subsets of $\hat{X}$ and $f\colon A \to B$ be a bijection. By Lemma \ref{thetasemiregularization} $A$ and $B$ are $\Theta$-discrete subsets of $X$. Hence $f$ can be extended to a homeomorphism $f\colon X \to X$ of $X$ which is also a homeomorphism of $\hat{X}$ since $f(\Int (\Cl (O))) = \Int (\Cl (f(O)))$ for any open subset $O$ of $X$. 
\end{proof}

In the light of the examples of Bing and Ritter, the following question looks interesting.

\begin{que}
\begin{itemize}
\item[(a)] Are there  examples of countable connected s$\Theta$-DH spaces, other than $B$ and $R$?
\item[(b)] Are there examples of countable connected Hausdorff sDH spaces?
\end{itemize} 

\end{que}

The following fact is relevant to this question.

\begin{pro}(\cite[Lemma 5]{MRV}) For every topological space $X$ with the topology $\tau$
we have $\hat{\hat{\tau}} =  \hat{\tau}$, i.e. the process of semi-regularizing of a 
topological space is an idempotent operation.
\end{pro} 

In particular, we have $\hat{R} = R.$

\section{Appendix}

The Bing space $B$ (\cite{B}) is the rational half plane $\{(x,y) \in \mathbb Q \times \mathbb Q\colon y \geq 0\}$ endowed with the topology $\tau$ consisting of the sets $U \subset \mathbb B$ such that for every  point $(a, b) \in U$ there exists $\epsilon > 0$ such that 

$\{(x, 0) \in B\colon | x-(a - b/\sqrt{3})| < \epsilon \} \cup
\{(x, 0) \in B\colon | x-(a + b/\sqrt{3})| < \epsilon \} \subset U. 
$

%\newpage

\vskip1cm

\noindent(V.A. Chatyrko)\\
Department of Mathematics, Linkoping University, 581 83 Linkoping, Sweden.\\
vitalij.tjatyrko@liu.se

\vskip0.3cm
\noindent(A. Karassev) \\
Department of Mathematics,
Nipissing University, North Bay, Canada\\
alexandk@nipissingu.ca

\end{document}